\documentclass[leqno,11pt]{amsart}

\usepackage{latexsym, esint, color}
\usepackage{amsmath,amssymb, amsthm}
\usepackage{graphicx}
\usepackage{caption}
\usepackage{subcaption}
\theoremstyle{plain}

\theoremstyle{definition}

\newcommand{\R}{\mathbb{R}}
\newcommand{\T}{\mathbb{T}}
\newcommand{\Z}{\mathbb{Z}}
\renewcommand{\S}{\mathbb{S}}
\renewcommand{\d}{\partial}
\renewcommand{\div}{{\rm div}\,}
\newcommand{\rot}{{\rm rot}\,}

\numberwithin{equation}{section}

\begin{document}

\begin{center}

{\Large Stabilizing the Long-time Behavior
of the Navier-Stokes Equations

\smallskip

\smallskip

and Damped Euler Systems
by Fast Oscillating Forces}

\bigskip

\bigskip

{\bf Jacek Cyranka$^{1,5}$, Piotr B Mucha$^{1}$, Edriss S Titi$^{2,3}$ and  Piotr Zgliczy\'nski$^4$ }

\smallskip

\smallskip

{\small 1. University of Warsaw, Institute of Applied Mathematics and Mechanics}

{\small ul. Banacha 2, 02-097 Warszawa, Poland}

\smallskip

{\small 2. Department of Mathematics, Texas A \& M University,}
{\small  3368 TAMU, College Station, TX 77843-3368, USA}

\smallskip

{\small 3. Department of Computer Science and Applied Mathematics,}
{\small  Weizmann Institute of Science, Rehovot 76100, Israel}

\smallskip

{\small  4. Institute of Computer Science, Jagiellonian University, \L ojasiewicza 6, 30-348 Krak\'ow, Poland}

\smallskip

{\small  5. Department of Mathematics, Rutgers, The State University of New Jersey, 110 Frelinghusen Rd, Piscataway, NJ  08854-8019, USA}

\smallskip

\smallskip

{Emails: cyranka@mimuw.edu.pl, \ p.mucha@mimuw.edu.pl, }

{ titi@math.tamu.edu, \ umzglicz@cyf-kr.edu.pl}

\end{center}
\date{December 16, 2015}

\vskip1cm

{\bf Abstract.}
 The paper studies the issue of stability of solutions to the Navier-Stokes and damped Euler systems in periodic boxes. We show that
 under action of fast oscillating-in-time external forces all two dimensional regular solutions converge to a time periodic flow.
 Unexpectedly, effects of stabilization can be also obtained for systems with stationary forces with large total momentum (average of the velocity).
 Thanks to the Galilean transformation and space boundary conditions, the stationary force changes into one with time oscillations.
 In the three dimensional case we show an analogical result for weak solutions to the Navier-Stokes equations.

 \vskip0,5cm

\section{Introduction}

In many analytical and computational studies of the forced Navier-Stokes or Euler equations, subject to periodic boundary conditions,
it is usually assumed that spatial average of the forcing term is zero. This in turn implies that the spatial average of the solution is invariant,
and for simplicity it is also taken to be zero. In this paper we investigate the long-time behavior of these systems when the spatial average of the initial 
velocity is taken to be large. By using the Galilean transformation of such  systems the problem is transformed into a similar system with fast time-oscillating
forcing term. Therefore, we investigate instead the long-time dynamics of the transformed Navier-Stokes equations and the damped Euler equations   under the action of
fast time-oscillating force. Specifically, we show that the
fast time-oscillating  forces have a stabilization effect;  and that the long-time dynamics   consists of a globally  attracting unique time-periodic solution.
This result is consistent with other results concerning the  investigation of the Navier-Stokes and Euler equation with  high oscillations. Fast rotation is, for
example,  a stabilizing mechanism of inviscid turbulent flows \cite{BMN,CG,GIM,MN}. Moreover, in the case of the forced two-dimensional Navier-Stokes equations it is observed that fast rotation is trivializing  the long-time dynamics, i.e., the  global attractor is a single stable steady state at the limit of large rotation rate \cite{W}.

Naturally, we distinguish in our proofs between the two-dimensional and the three-di\-men\-sio\-nal cases. For the two-dimensional models  we show that the long-time
dynamics of regular solutions is trivial; specifically,  that all solutions tend to a unique time-periodic solution
generated by the fast time-oscillating forcing term. For the three-dimensional  case we  only consider the Navier-Stokes equations,  and we show that all Leray-Hopf weak solutions
converge, as time tends to infinity, to a unique time periodic solution generated by the fast time-oscillating forcing term.
Notably, we do not impose any  restriction on the magnitude of spatial norms of the forcing term, and we only assume that the size spatial average of the initial data (or equivalently the rate of oscillation in the forcing term) is large enough. In the last section of this paper we also provide some numerical results supporting our qualitative theoretical results.

The paper is organized as follows. In the next section we motivate our study and introduce the relevant models.
In section 3 we treat the two-dimensional  models, and in section 4 we consider the three-dimensional  Navier-Stokes system. Numerical results are reported in section 5. We also provide an Appendix section in which we sketch the construction of time periodic solutions.

\section{Settings  and motivation}

In this paper we consider the Navier-Stokes and the damped Euler   systems of equations
\begin{equation}\label{EeS}
 \begin{array}{l}
  v_t + v \cdot \nabla v -\epsilon \Delta v + \alpha [v] + \nabla p= F,\\[6pt]
\div v =0 ,\\
 \end{array}
\end{equation}
subject to periodic boundary conditions in the  $n$-dimensional torus $\T^n$, for $n=2,3$, and with divergence-free initial datum $v_0$. Here we denote by
\begin{equation}\label{i2a}
 [v]=v -\frac{1}{|\T^n|} \int_{\T^n} v dx,
\end{equation}
and assume that
\begin{equation}\label{i2}
 \int_{\T^n} F(x,t) dx =0, \mbox{ \ \ for all \ }  t >0.
\end{equation}
Assumptions \eqref{i2a} and \eqref{i2} imply  the conservation of the total momentum of the flow, i.e.,
\begin{equation}\label{i3}
 \int_{\T^n} v(x,t) dx = \int_{\T^n} v_0(x)dx.
\end{equation}
Throughout this work  we assume that $\max\{ \alpha, \epsilon\} >0$.

In this section we consider two special cases of system \eqref{EeS}:

\smallskip

\noindent {\bf Case (A):}\  In the first case we consider a very fast (i.e., $\Omega$ in \eqref{i4} below is very large)  constant background  flow motion in the  $x_1-$direction, given by the initial data, i.e.,
\begin{equation}\label{i4}
\frac{1}{|\T^n|} \int_{\T^n} v_0 dx =  \Omega \hat{e}_1.
\end{equation}
In addition, we assume that the external force in \eqref{EeS}  is time independent, i.e.,  $F(x,t)=F(x)$;  and that the
 Fourier  coefficients $\hat{F}_{k}=0$, for $k=(0,k_2)$ in the 2d case, and  for $k=(0,k_2,k_3)$ in the 3d case. As it will become clear later, this assumption implies  that the forcing term does not resonate with the background constant flow given in \eqref{i4}.

The choice of the direction $x_1$ for the background flow is not important, but it simplifies the presentation.
The above assumptions allow  us in this case to make the following change of the variables, using the Galilean transformation,
\begin{equation}\label{i5}
 x \rightarrow x + \Omega t \hat{e}_1  \mbox{ \ \ and \  \ } v \to v - \Omega \hat{e}_1 \mbox{ \ with \ } \Omega \in \R_+.
\end{equation}
Consequently,  we arrive to the following equivalent system to  (\ref{EeS})


\begin{equation}\label{EeS2}
 \begin{array}{l}
  v_t + v \cdot \nabla v -\epsilon \Delta v + \alpha v + \nabla p= F(x_1+\Omega t,x'),\\[6pt]
\div v =0,
 \end{array}
\end{equation}
with initial datum $v_0$ such that
$ \int_{\T^n} v_0 dx =0.$ Here  $x'=x_2$ for the 2d case, and $x'=(x_2,x_3)$ for the 3d case. For large values of $\Omega$,  the speed of the background flow, transformation (\ref{i5}) yields a new system, \eqref{EeS2}, that governs the perturbation about the background flow, with fast oscillating in time forcing term, with period $T_{per} =\frac{2\pi}{\Omega}$, but with zero total momentum.

\smallskip

\noindent{\bf Case (B):} \  In the  second case we consider system \eqref{EeS} with a  special type of fast oscillating forcing term.
Specifically,  we consider fast oscillating force of the form
\begin{equation}\label{a1}
 F(x,t)= f(x) \sin \Omega t.
\end{equation}
We also assume that
$\displaystyle \int_{\T^n} f(x) dx =\int_{\T^n} v_0(x) dx =0.$ Here $[v]=v$.

\smallskip

The key observation in both, case (A) and case (B),  is that we consider systems with fast oscillating forcing terms in (\ref{EeS2}) and (\ref{a1}), respectively. The main purpose of this study is to take advantage of these fast oscillating forcing terms  to stabilize the long-time behavior of the solutions of the corresponding systems. Consequently, our analysis will concentrate on the limit, as $\Omega \to \infty$. Note that the force (\ref{a1}) is, roughly speaking, a particular case of the force considered in  (\ref{EeS2}),  since $F(\cdot)$ is spatially periodic.

\smallskip

In both cases, (A) and (B),  we require the force $F$ to be sufficiently smooth and we do not restrict the magnitude of  its spatial norms.
In  particular, the $L_2(\T^n)$ and $H^{-1}(\T^n)$ type norms
can be arbitrary large and fixed. The same we assume about the size of the initial data.

\smallskip

Next, we provide a rough description of  the main results presented in this article.

\smallskip

\noindent {\bf Results for the two-dimensional  case:}

\smallskip

{\it 1. Periodic in time solutions.} Let $\delta>0$ be sufficiently small, and let $F$ be sufficiently smooth, fulfilling the forms in case (A) or case (B).
Let $\epsilon \geq 0$ and $\alpha \geq 0$ with $\max\{\epsilon,\alpha\}>0$. Then there exists
a  periodic  in time solution, $v_{per}$, to (\ref{EeS}) such that
\begin{equation}\label{a2}
 \|\nabla v_{per}\|_{L_\infty(\T^2 \times T_{per}\S^1)} \leq \delta,
\end{equation}
provided $\Omega$ is sufficiently large. The period is  $T_{per}=2\pi / \Omega$. Theorem 1 from Section 3.

\smallskip

{\it 2. Global stability and uniqueness of the periodic solution}. Let $v$ be a solution to (\ref{EeS}), with arbitrary, divergence-free, initial datum $v_0\in L_2$. Then for every  $\Omega$  large enough, such that the above statement  is valid, we have
\begin{equation}\label{a3}
 \|v(t)-v_{per}(t)\|_{L_2(\T^2)} \to 0, \mbox{ as } t \to \infty,
\end{equation}
where $v_{per}$ is as above. In particular, and under the above assumptions, it follows from  \eqref{a3}  that $v_{per}$ is unique.
Theorem 2 from Section 3.

Here we shall mention about a current result from \cite{CZ} concerning dissipative PDEs which intersects with the case A for some class of
external forces. Note however that methods used there are significantly different from ones we apply in the present paper.

\smallskip

\noindent
{\bf Results for the three-dimensional case:}

\smallskip

We also investigate the long-time behavior of the three-dimensional Navier-Stokes equations, for large values of $\Omega$.   However,  due to our inability to prove global existence of strong solutions, or the uniqueness of weak solutions for the 3d Navier-Stokes equations, we will focus on the long-time behavior of the Leray-Hopf class of weak solutions
 of the 3d Navier-Stokes system. Thus, by taking $\alpha =0$ and $\epsilon =\nu >0$ in (\ref{EeS}),
we consider the following 3d Navier-Stokes system:
\begin{equation}\label{NSE}
 \begin{array}{l}
  v_t + v \cdot \nabla v -\nu \Delta v  + \nabla p= F,\\[6pt]
\div v =0, \\
 \end{array}
\end{equation}
subject to periodic boundary conditions on the torus  $\T^3$,  and with  the divergence-free initial datum $v_0\in L_2(\T^3)$ that has zero spatial average on $\T^3$. Moreover, the  forcing term in \eqref{NSE}, $F$, is assumed  to satisfy the  conditions in either case (A) or case (B), above.

\smallskip

{\it 3. Periodic in time solutions.} Let $\delta>0$ be sufficiently small. Let $F$ be sufficiently smooth satisfying either case (A) or case (B).   Then there exists
a  periodic  in time solution, $v_{per}$,  to (\ref{NSE}) such that
\begin{equation}\label{a4}
 \| v_{per}(t) \|_{L_\infty(\T^3 \times T_{per}\S^1 )} \leq \delta,
\end{equation}
provided $\Omega$ is sufficiently large. Theorem 3 from Section 4.

\smallskip

{\it 4. Global stability and uniqueness of the periodic solution}. Let $v$ be a
Leray-Hopf weak  solution to (\ref{NSE}), with arbitrary divergence-free initial datum $v_0\in L_2(\T^3)$. Then, for  every $\Omega$  sufficiently large, we have
\begin{equation}\label{a5}
 \|v(t)-v_{per}(t)\|_{L_2(\T^3)} \to 0, \mbox{ as } t \to \infty,
\end{equation}
where $v_{per}$ is as above. In particular, and under the above assumptions, it follows from  \eqref{a5}  that $v_{per}$ is unique.
Theorem 4 from Section 4.

\smallskip

\noindent{\bf  Numerical results for the two-dimensional case:}

We conclude the paper with numerical tests  illustrating our analytical results for certain class of flows.
Specifically, we consider the 2D Kolmogorov flows in the flat torus $\mathbb{T}^2_\beta =[0,2\pi]\times[0,2\pi \beta]$, for $\beta >0$; which is  the  2D NS equations, subject to periodic boundary condition, forced by an eigen-function of the Stokes operator.   This problem is on the one hand simple, but on the other hand  is dynamically rich enough   to illustrate the phenomena at hand.

We consider the vorticity formulation of the 2D version of \eqref{EeS2} with the specific  forcing term $F=(-\lambda\beta \sin(\frac{y}{\beta} + \Omega t), 0)$ which yields

\begin{equation}
\label{eq:omegaExpOld}
\omega_t+v\cdot\nabla\omega-\epsilon\Delta\omega+\alpha\omega=
    \lambda\cos(\frac{y}{\beta} + \Omega t ).
\end{equation}

 Observe that $\omega(z,t)$ is a solution of \eqref{eq:omegaExpOld}, for  $z=(x,y)\in\mathbb{T}^2_\beta$, if and only if $\tilde{\omega}(z,t)=\omega(z+\beta\Omega t \hat{e}_2,t)$ is a solution to the evolution equation
\begin{equation}
\label{vortSystem}
\tilde{\omega}_t+(\tilde{v} + \beta\Omega \hat{e}_2)\cdot\nabla\tilde{\omega}-\epsilon\Delta\tilde{\omega}+\alpha\tilde{\omega}=
    \lambda\cos(\frac{y}{\beta}).
\end{equation}
In particular, $\omega^*(z)$ is a stationary solution of \eqref{vortSystem} if and only if $\omega(z,t) = \omega^*(z+\beta\Omega t \hat{e}_2)$, for
$z=(x,y)\in\mathbb{T}^2_\beta$, is a   time periodic solution  to \eqref{eq:omegaExpOld}. Moreover,  $\omega^*(z)$ is globally stable for the dynamics of \eqref{vortSystem} if and only if $\omega^*(z+\beta\Omega t \hat{e}_2)$ is global stable time periodic solution of \eqref{eq:omegaExpOld}.

\begin{itemize}
\item Based on the above observation we  present a bifurcation analysis of the stationary solutions to \eqref{vortSystem}. In particular,
 investigate the bifurcation diagram  of  \eqref{vortSystem}  for large values of $\Omega$. Moreover, we show that the range of $\lambda$'s for which   system \eqref{vortSystem} admits a unique stationary  solution increases proportionally to $\Omega$.

\item
We investigate the stationary solutions of   \eqref{vortSystem} rather than direct numerical integration   in  order to avoid working with rapidly
oscillating in time functions.

\item We also investigate the rate of convergence to the globally stable solution of \eqref{vortSystem}. The purpose of this study is to    provide an evidence  of
the exponential convergence rate, which we show in Theorem~2.
\end{itemize}

All of the numerical results were derived using a finite dimensional Galerkin approximations, we argue that the dimensions we used are sufficient.

\smallskip

\smallskip

\noindent
{\bf An illustrative linear toy model with friction --  Newton's second law}

\smallskip

To illustrate the stabilization mechanism, due to the fast oscillations in the forcing term, we focus here on the following simple  linear equation with friction/damping/drag term:
\begin{equation}\label{s1}
 w_t + \alpha w = F(x,\Omega t),
\end{equation}
 in the torus $\T^n$. Here $F$ is time periodic, with period $T_{per} = \frac{2\pi}{\Omega}$.
First, we observe  that the solution to  system (\ref{s1}) does not involve dynamically  the spatial variable, $x$, so the solution will treat $x$  as
a parameter (label), i.e., we have a  parameterized system  of simple
ODEs.

We assume that the forcing term $F(x,\Omega t)$, in \eqref{s1}, enjoys specific structure, namely,  there exists a smooth function $g(x,\Omega t)$,
periodic in space and time,  such that
\begin{equation}\label{s2}
 \d_t g(x, t) =  F(x, t),  \mbox{\ consequently \ \ \ }\frac{1}{\Omega} \d_t g(x,\Omega t) =  F(x,\Omega t).
\end{equation}
We have two prototypical examples in mind of the forcing terms, $F(x,\Omega t)$, satisfying  the above structure. Specifically, let   $f(x)$, for $x\in \T^n$, be a smooth spatially periodic function. We consider again the cases:
\begin{equation}\label{s3}
 (A) \quad g(x,\Omega t) = \frac{1}{\Omega}D^{-1}_1f(z)|_{z=x+\Omega t \hat{e}_1},
 \qquad   (B) \quad g(x,\Omega t) = - \frac{1}{\Omega } f(x) \cos \Omega t.
\end{equation}
Here we put $D^{-1}_1f(z)$ as the primitive function of $f$ with respect to the first variable, i.e.,
we have $\frac{1}{\Omega} \d_t D^{-1}_1f(x+\Omega t \hat{e}_1)=f(x+\Omega t \hat{e_1})$.
Put it in other words we define
\begin{equation}\label{s3a}
g(x,\Omega t) =\sum_{k \in \Z^n} \frac{\hat{f}_k}{ik_1} e^{ik\cdot x} e^{i\Omega t k_1}
\end{equation}
 via Fourier series. Here we see that the assumption
$\hat{f}_k=0$,  for  $k=(0,k')$, is required to justify the above form of $g$.

As a result of the previous assumptions on  $f$ we have in both cases that
\begin{equation}\label{s4}
 \|g(\cdot,\Omega t)\|_{C^2(\T^n)} \leq \frac{C}{\Omega},
\end{equation}
where $C$ depends on the spatial norms of $f$, but is independent of $\Omega$. Notice that
only the time derivatives of $g$ will add  multiplication by factors of $\Omega$.
 Set
\begin{equation}\label{s5}
 w(x,t)= W(x,t) + g(x,\Omega t),
\mbox{ \ \ then $W$ satisfies \ \ }
 W_t + \alpha W = \alpha g.
\end{equation}
Therefore,  from (\ref{s5}), by (\ref{s4}) we have
\begin{equation}\label{s6}
 W(x,t)= \exp\{ -\alpha t\} W_0(x) + \int_0^t \exp\{-\alpha(t-s)\} \alpha g(x,\Omega s) ds \sim e^{-\alpha t} + \frac{1}{\Omega}.
\end{equation}
Thus, the solution to (\ref{s1}) satisfies
\begin{equation}
\| w (t)\|_{L_\infty(\T^n)} \sim e^{-\alpha t} + \frac{1}{\Omega}.
\end{equation}
 Since the problem is linear, the above structure  holds for arbitrary positive $\alpha$ and $\Omega$.

Next, let us consider the time periodic solutions to   (\ref{s5}):
\begin{equation}\label{s5a}
 W_{per,t} +\alpha W_{per} =\alpha g.
\end{equation}
The construction of periodic solutions to the \eqref{s5a} can be done explicitly through the Fourier series in time, on the time periodic interval  $T_{per}\S$. Using the energy methods we immediately obtain
the following bound
\begin{equation}
 \|W_{per}\|_{L_\infty(T_{per}\S^1)} \leq \|g\|_{L_\infty(T_{per}\S^1)}.
\end{equation}


%
Comparing the solutions to (\ref{s5}) and to (\ref{s5a}) we obtain the trivial identity
\begin{equation}
 \d_t(W - W_{per}) + \alpha (W - W_{per}) =0.
\end{equation}
The above identity implies
\begin{equation}
 |W (x,t)- W_{per}(x,t)| = |W_0(x)-W_{per}(x,0)| e^{-\alpha t}\sim  e^{-\alpha t}.
\end{equation}

Summing up we obtain the following:

\smallskip

{\it {\bf Proposition 1.} Let $\alpha >0$, and $g$ has one of the forms in (\ref{s3}), then every solution to (\ref{s1}) admits
the following structure
\begin{equation}
 \|w(\cdot,t)\|_{L_\infty(\T^n)} \lesssim e^{-\alpha t}\|w_0\|_{L_\infty(\T^n)}
 +\frac{1}{\Omega}\|f\|_{L_\infty(\T^n)}.
\end{equation}
and
\begin{equation}
  \|w(\cdot,t)-w_{per}(\cdot,t)\|_{L_\infty(\T^n)} \leq e^{-\alpha t}\|w_0-w_{per}(\cdot,0)\|_{L_\infty(\T^n)},
\end{equation}
where $w_{per}=W_{per}-g$.
}

\section{The two-dimensional case}

System (\ref{EeS}), that we consider here, is a modification of the Navier-Stokes and  Euler equations, by basically adding a linear friction/damping/drag
force with coefficient $\alpha \ge 0$. We require in addition that $\max\{ \alpha,\epsilon\} >0$,  thus, (\ref{EeS}) is a dissipative form of the  Euler system.
 We use the special proprieties which are valid in the 2d case.
Namely, we analyze  system (\ref{EeS})
in the vorticity formulation which takes the form
\begin{equation}\label{b2}
 \omega_t + v \cdot \nabla \omega - \epsilon \Delta \omega +\alpha \omega = \rot F,
\end{equation}
where
\begin{equation}\label{b1}
 \omega = {\rm rot}\, v.
\end{equation}
Our result concerning  system (\ref{b2})-(\ref{b1}) is the following

\smallskip

{\it {\bf Theorem 1.} Let $\alpha \ge 0$, $\epsilon \ge 0$, with $\max\{\alpha,\epsilon\} >0$, and  let $F$ be sufficiently smooth force
of form (A) or (B).
In addition, for case (A) let us assume that there exists a scalar function $g(x,\Omega t)$ such that
$$
\frac{1}{\Omega} \d_t  g(x,\Omega t)=\rot F(x_1+\Omega t,x_2) \mbox{ \ \ and \ \ }
\sup_{s} \|g(\cdot,s)\|_{C^3(\T^2)} \leq G,
$$
where $G$ is independent of $\Omega$.
Choose $\delta$ small enough such that
$$
0< \delta \leq \frac 14 (\alpha+ a \epsilon),
$$
where $a$ is an absolute constant.
Then there exists a regular periodic in time solution to problem (\ref{EeS}) such that
\begin{equation}\label{b6}
 \|v\|_{L_\infty(T_{per}\S^1;W^1_\infty(\T^2))} \leq \delta,
\end{equation}
provided $\Omega$  is sufficiently large, depending on $\alpha,\epsilon,G$ and $\delta$, with $T_{per}={2\pi}/{\Omega}$.
}

\smallskip

{\bf Proof.} The construction of time periodic solutions is based on the domain $\T^2 \times T_{per}\S^1$. We consider
system (\ref{EeS}) in the form of (\ref{b2}). Let $\bar v$ be a given smooth enough time periodic velocity field satisfying:
\begin{equation}\label{b9}
 \div \bar v=0 \mbox{ \ with \ } \|\bar v\|_{L_\infty(T_{per}\S^1;W^1_\infty(\T^2))} \leq \delta.
\end{equation}
We then look for a time periodic  vorticity $\omega$ that solves the following "linearized" transport, by the velocity field $\overline{v}$, version of   problem (\ref{b2}):
\begin{equation}\label{b8}
 \omega_t + \bar v \cdot \nabla \omega  - \epsilon \Delta \omega +\alpha \omega = \rot F.
\end{equation}
Then we construct the time periodic  velocity field, $v$,  corresponding to the vorticity $\omega$ such that
\begin{equation}\label{b10}
 \rot  v = \omega, \qquad \div  v =0, \qquad \int_{\T^2} v dx=0.
\end{equation}
 We show that the above procedure defines a map
$$
\mathcal{T}:L_\infty(T_{per}\S^1;W^1_\infty(\T^2)\cap \{\div v=0\}) \to L_\infty(T_{per}\S^1;W^1_\infty(\T^2)\cap \{\div v=0\}),
$$
 such that $\mathcal{T}(\bar v)=v$.  And we look for a fixed point of this map, which will in turn define a time periodic solution to (\ref{b2}).  Indeed, we  show that $\mathcal{T}$ maps the set
$$
\Xi=\{\div v=0 \mbox{ with } \|v\|_{L_\infty(T_{per}\S^1;W^1_\infty(\T^2))} \leq \delta\}
$$
into itself and  that  $\mathcal{T}$ is a compact map. Then the assumptions of the Schauder  fixed point theorem are
fulfilled, which will imply the existence of a fixed point
of map $\mathcal{T}$.
Given $\bar v$ satisfying (\ref{b9}), the existence of time periodic solution to the  linearized system (\ref{b8}) is not difficult, and it can be proved easily by the Galerkin method -- see the Appendix.
Next, we establish the required estimates.

\smallskip

We split our proof into two cases, distinguishing influences of $\epsilon \Delta w $ and  $\alpha w$.


\bigskip
\noindent {\it The dominant-damping case} is when  $\alpha \geq a\epsilon$, for some positive absolute constant $a$ to be specified later. This is the case when  the damping $\alpha>0$ is dominating the  viscosity $\epsilon$ which  is very small or maybe even equal to zero. In the latter case we essentially  have the damped Euler system. Recalling (\ref{s3}) we set
\begin{equation}\label{b11}
 W=\omega -\frac{1}{\Omega} g(x, \Omega t).
\end{equation}
 For case (A),  the function $g$ is given by (\ref{s3a}); and for  case (B) we take, as before,  $g(x,\Omega t )=- \cos \Omega t \cdot f(x)$.
Then $\d_t \frac{1}{\Omega} g(x,\Omega t) = \rot F$.
Then $W$ fulfills
\begin{equation}\label{b12}
 W_t + \bar v \cdot \nabla W -\epsilon \Delta W +\alpha W = \bar v \cdot \nabla \frac{1}{\Omega} g  -  \epsilon \Delta
\frac{1}{\Omega} g + \alpha \frac{1}{\Omega} g.
\end{equation}
Multiplying  (\ref{b12})  by $|W|^{p-2}W$  and integrating over $\T^2$ we obtain

\begin{multline}\label{b13}
 \frac{1}{p} \frac{d}{dt} \|W\|^p_{L_p(\T^2)} + \epsilon \int_{\T^2} (p-1) |\nabla W|^2 |W|^{p-2} dx
 + \alpha \|W\|^p_{L_p(\T^2)} dx \leq \\
\frac{1}{\Omega} \int_{\T^2} \big( |\bar v| |\nabla g| + \epsilon |\Delta g| + \alpha |g|\big) |W|^{p-1} dx.
\end{multline}
Integrating with respect to time, over $T_{per}\S^1$, and using the periodicity in time,  we get
\begin{multline}\label{b14}
 \epsilon \int_{T_{per}\S^1} \int_{\T^2} (p-1) |\nabla W|^2 |W|^{p-2} dxdt  + \alpha \|W\|^p_{L_p(\T^2\times T_{per}\S^1)}
 \leq \\
C (\frac{G}{\Omega} \| \bar v\|_{L_p(T_{per}\S^1; L_\infty(\T^2))} \|W\|_{L_p(T_{per}\S^1; L_p(\T^2))}^{(p-1)/p}
+ \frac{G}{\Omega}(\epsilon + \alpha) T_{per}^{1/p} \|W\|_{L_p(T_{per}\S^1; L_p(\T^2))}^{(p-1)/p} ).
\end{multline}
In particular, we have


\begin{equation}\label{b15}
 \alpha \|W\|_{L_p(T_{per}\S^1; L_p(\T^2))} \leq
 C\frac{G}{\Omega} T^{1/p}_{per}(\| \bar v\|_{L_\infty(T_{per}\S^1; L_\infty(\T^2))} + \alpha).
\end{equation}
However this regularity is not enough, one more spatial derivative is required, so we differentiate (\ref{b12}) with respect to $x_i \in \{x_1,x_2\}$ getting
\begin{multline}\label{b16}
 W_{x_i,t} + \bar v \cdot \nabla W_{x_i} -\epsilon \Delta W_{x_i} +\alpha W_{x_i} = -\bar v_{x_i} \cdot \nabla W
\\
+ \bar v_{x_i} \cdot \nabla \frac{1}{\Omega} g + \bar v \cdot \nabla \frac{1}{\Omega} g_{x_i}-  \epsilon \Delta
\frac{1}{\Omega} g_{x_i} + \alpha \frac{1}{\Omega} g_{x_i}.
\end{multline}
Multiplying  (\ref{b16}) by $|W_{x_i}|^{p-2}W_{x_i}$, integrating over $\T^2\times T_{per}\S^1$, and using the periodicity in time, we get
\begin{multline}\label{b17}
 \sum_{i=1}^2 \big[(p-1)\epsilon \int_{T_{per}\S^1}\int_{\T^2} |\nabla W_{x_i}|^2 |W_{x_i}|^{p-2} dxdt + \alpha \int_{T_{per}\S^1} \int_{\T^2} |W_{x_i}|^p dxdt \big] \\
\leq
\|\nabla \bar v\|_{L_\infty} \|\nabla W\|^p_{L_p(T_{per}\S^1; L_p(\T^2))} + \frac{G}{\Omega} T_{per}^{1/p} ( \| \bar v\|_{L_\infty(T_{per}\S^1; W^1_\infty(\T^2))}
 +\alpha) \|\nabla W\|^{(p-1)/p}_{L_p(T_{per}\S^1; L_p(\T^2))}.
\end{multline}
Assuming, as required in (\ref{b6}),   that  $\|\nabla \bar v\|_{L_\infty} \leq \delta$, and observing that  in this case we have $ \delta \leq \frac 12 \alpha$, we conclude that
\begin{equation}\label{b18}
 \alpha \|\nabla W\|_{L_p(T_{per}\S^1; L_p(\T^2))} \leq
C \frac{G}{\Omega} T^{1/p}_{per}( \| \bar v \|_{L_\infty(T_{per}\S^1; W^1_\infty(\T^2))}  + \alpha) \leq C\frac{\alpha G}{\Omega} T^{1/p}_{per}.
\end{equation}
Substituting estimates (\ref{b18}) and (\ref{b15}) into  equation (\ref{b12}) we find that
\begin{equation}\label{b19}
 W_t - \epsilon \Delta W \in L_p(T_{per}\S^1; L_p(\T^2)).
\end{equation}
Based on the classical result for the heat equation \cite{A,LSU}
of the maximal regularity estimates for the $L_p$ spaces, we obtain information with no dependence  from  $\epsilon$. Hence
\begin{equation}\label{b20}
 \|W_t\|_{L_p(T_{per}\S^1;L_p(\T^2))} \leq
 C\frac{G}{\Omega} T^{1/p}_{per} ( \| \bar v\|_{L_\infty(T_{per}\S^1; W^1_\infty(\T^2))} +  \alpha)
 \leq C\frac{\alpha G}{\Omega} T^{1/p}_{per}.
\end{equation}
 Concerning estimates (\ref{b12})-(\ref{b20}), we observe that the term $-\epsilon \Delta W$, for $\epsilon >0$, has the right sign. Moreover, it also has the right  good sign even in the maximal regularity (\ref{b19}).

\bigskip

\noindent{\it The dominant-viscosity case} when $a\epsilon \geq \alpha$, with possibly  $\alpha=0$. This case allows us to take full advantage of the
parabolicity of (\ref{b12}), which reads
\begin{equation}\label{b40}
 W_t  -\epsilon \Delta W +\alpha W = -\bar v \cdot \nabla W+ \bar v \cdot \nabla \frac{1}{\Omega} g  -  \epsilon \Delta
\frac{1}{\Omega} g + \alpha \frac{1}{\Omega} g.
\end{equation}
The maximal regularity estimate for the heat equation \cite{LSU} or more direct \cite{Mu2} implies that

\begin{equation}\label{b41}
 \|W_t,\epsilon \nabla^2 W, \alpha W\|_{L_p(\T^2 \times T_{per}\S^1)} \leq C_p\| \bar v \cdot \nabla W, \bar v \cdot \nabla \frac{1}{\Omega} g,    \epsilon \Delta
\frac{1}{\Omega} g, \alpha \frac{1}{\Omega} g\|_{L_p(\T^2 \times T_{per}\S^1)} ,
\end{equation}
where the constant $C_p$ depends only on $p$, there is no dependence on $T_{per}$, since we consider only homogeneous norms in (\ref{b41}).
Observe that
\begin{equation*}
 C_p\| \bar v \cdot \nabla W, \bar v \cdot \nabla \frac{1}{\Omega} g,    \epsilon \Delta
\frac{1}{\Omega} g, \alpha \frac{1}{\Omega} g\|_{L_p(\T^2 \times T_{per}\S^1)} \leq C_p\delta \|\nabla W\|_{L_p\T^2 \times T_{per}\S^1)}
+ C\frac{\epsilon G}{\Omega} T^{1/p}_{per},
\end{equation*}
since $\bar v \in \Xi$
we have $\|\bar v\|_{L_\infty} \leq \delta$.
Furthermore, since in this case we have $\delta \leq \frac 12 a \epsilon$ then the first term above can be absorbed by  the left-hand side  of (\ref{b41}),
 thanks to the facts $\int_{\T^2} W dx =0$ and $\|\nabla W\|_{L_p(\T^2)} \leq C\|\nabla^2 W\|_{L_p(\T^2)}$.
Hence we conclude the following bound


\begin{equation}\label{b42}
 \|W_t,\epsilon \nabla^2 W, \alpha W\|_{L_p(\T^2 \times T_{per}\S^1)} \leq C\frac{\epsilon G}{\Omega} T^{1/p}_{per},
\end{equation}
for the case $a\epsilon \geq \alpha$, establishing  the analogue  of (\ref{b18}) and (\ref{b20}) for this case.

\smallskip

Now we return to studying properties of the map $\mathcal{T}$ treated for both cases. Before we establish the $L_\infty$-bound for $\nabla v$, a comment is in order.
The key problem is the length of $T_{per}$. In general the constant in the Sobolev imbeddings
may highly depend on $T_{per}$ in a bad way. Hence a solution, which here seems to be most natural, is to consider
 $\bar v$ over several  periods of time. Here we think about $\S^* \sim [\frac{1}{T_{per}}]T_{per}\S^1$, where $[t]$ denotes the integer part of $t$.
Then $\S^*$ is close to $\S^1$. The functions are defined over the domain $\T^2 \times \S^*$, so the problems with thinness of domain will
disappear.

 For this purpose we set  $v=v_2-v_1$,  where   $v_2$ is
given as a solution to the following problem
\begin{equation}\label{bb1}
 \rot v_2 = W, \qquad \div v_2=0
\end{equation}
and  $v_1$ is given by
\begin{equation}\label{bb2}
 \rot v_1 = \frac{1}{\Omega} g(x,\Omega t), \qquad \div v_1=0.
\end{equation}
The functions are considered on time interval $\S^*$, since we assumed that $\Omega$ is large, hence $T_{per} <<1$.
Then from (\ref{b18}) and (\ref{b20}), together with (\ref{b15}), and from (\ref{b42}), we get


\begin{equation}\label{bb4}
 \|\nabla v_2\|_{W^{1,1}_p(\T^2 \times \S^*)} \lesssim
\frac{1}{T_{per}^{1/p}} \|\rot v_2\|_{W^{1,1}_p(\T^2 \times T_{per}\S^1)} \leq C(1+\alpha)\frac{G}{\Omega} \lesssim \frac{G}{\Omega}.
\end{equation}

So the Sobolev imbedding $W^{1,1}_p (\T^2 \times \S^*) \subset L_\infty(\T^2 \times \S^*)$  gives
\begin{equation}\label{bb5}
 \|\nabla v_2 \|_{L_\infty(0,1;L_\infty(\T^2))} \leq C\frac{G}{\Omega},
\end{equation}
where $C$ is independent of $T_{per}$.
Thus,
\begin{equation}
 \|\nabla  v \|_{L_\infty(\T^2 \times \S^*)} \leq
 \|\nabla v_2 \|_{L_\infty(\T^2\times \S^*)} + \|\nabla v_1 \|_{L_\infty(\T^2 \times \S^*)}\leq C\frac{G}{\Omega} \leq \delta,
\end{equation}
provided $\Omega$ large enough.

We showed that $\mathcal{T}$ maps the set $\Xi$ into itself, and the imbedding (for the case $\alpha \geq a \epsilon$)
\begin{equation}\label{b21}
 W^{1,1}_p(\T^2 \times T_{per}\S^1) \subset L_\infty(\T^2 \times T_{per}\S^1)
\end{equation}
 for $p>3$ is compact -- (\ref{b18}) and (\ref{b20}). The space $W^{1,1}_p(\T^2 \times T_{per}\S^1)$ is defined as a set of functions $f$ such that
$\nabla_x f\in L_p(\T^2 \times T_{per}\S^1)$ and $\d_t f\in L_p(\T^2 \times T_{per}\S^1)$. The case $a\epsilon \geq \alpha$ is simpler.

By  the Schauder fixed theorem we obtain existence of at least one fixed point
 of the map $\mathcal{T}$ fulfilling (\ref{b6}).
This implies existence of time periodic solutions to the nonlinear system (\ref{b2})-(\ref{b1}) satisfying bound (\ref{b6}). Theorem 1 is proved.

\smallskip

\noindent
{\bf Global stability of the time periodic solutions}

\smallskip


{\it {\bf Theorem 2.} Let $p>3$, let $v(t)$ be a solution to (\ref{b2}-\ref{b1})  corresponding to the initial data $v_0 \in W^2_p(\T^2)$ and $F$ fulfills the assumptions
 from Theorem 1, then
\begin{equation}\label{c1}
 \|v(t)-v_{per}(t)\|_{L_2} \to Ce^{-\frac{(\alpha +\epsilon)  t}{2}}, \mbox{ as } t \to \infty,
\end{equation}
where $v_{per}$ is the time periodic solution established by Theorem 1. Moreover $v_{per}$ must be unique.
}

\smallskip

{\bf Proof.} Note that for smooth enough initial datum $v_0$ we have the  global in time existence of solutions to system (\ref{EeS}).  Consider the difference
$\delta v(t)= v(t) - v_{per}(t),$ it fulfills the following system
\begin{equation}\label{c2}
 \begin{array}{l}
  \delta v_t + v \cdot \nabla \delta v -\epsilon \Delta \delta v + \alpha \delta v + \nabla \delta p= -\delta v \cdot \nabla v_{per} ,\\[6pt]
\div \delta v =0, \\[6pt]
\delta v|_{t=0}=v_0 - v_{per}(\cdot,0).
 \end{array}
\end{equation}
Multiplying  (\ref{c2}) by $\delta v$ and integrating over $\T^2$ yields
\begin{equation}\label{c3}
 \frac{1}{2}\frac{d}{dt} \int_{\T^2} (\delta v)^2 dx+ \int_{\T^2} \big(\epsilon |\nabla \delta v|^2 + \alpha (\delta v)^2 \big) dx \\
\leq \|\nabla v_{per}\|_{L_\infty}  \int_{\T^2} (\delta v)^2dx.
\end{equation}
Applying the Poincar\'e inequality we get
\begin{equation}\label{c3a}
 \frac{1}{2}\frac{d}{dt} \int_{\T^2} (\delta v)^2 dx+ (\epsilon +\alpha) \int_{\T^2}  (\delta v)^2  dx
\leq \delta \int_{\T^2} (\delta v)^2dx.
\end{equation}
In our setting the constant from the Poincar'e inequality (in the $L_2$ spaces) is equal 1.
Our choice of $\delta$ guaranteed that $\delta \leq \frac{1}{2}(\epsilon +\alpha)$,
so  we get
\begin{equation}\label{c4}
 \frac{d}{dt} \int_{\T^2} (\delta v)^2 dx+ (\epsilon +\alpha) \int_{\T^2}  (\delta v)^2 dx
\leq 0.
\end{equation}
We immediately conclude (\ref{c1}).
In particular (\ref{c1})  shows that constructed  time periodic solution  by Theorem 1 is unique.

\section{The 3d case}

{\it {\bf Theorem 3.} Let $\nu > 0$, and $F$ be sufficiently smooth divergence-free vector field of form (A) or (B).  Suppose that  for case (A) there   exists a divergence-free vector field
 $H$ satisfying
\begin{equation}\label{d1}
 \d_t \frac{1}{\Omega} H(x,\Omega t) = F(x_1+\Omega t,x') \mbox{ \ and \ } \sup_t \|H(\cdot, \Omega t)\|_{C^1(\T^3)} \leq G.
\end{equation}
There exists $\Omega_0(\nu,G)>0$ such that
 there exists a time periodic solution, with period $T_{per}=2\pi/\Omega$ satifying
\begin{equation}\label{d2}
 \|v_{per} \|_{L_\infty(\T^3 \times T_{per}\S^1)} \leq \frac{{C^*}(1+\nu)G}{\Omega},
\end{equation}
provided $\Omega \geq \Omega_0$, where ${C^*}$ depends on $\nu$, only.}

\smallskip

{\bf Proof.} Let $\delta=\frac{{C^*}(1+\nu)G}{\Omega}$.
Introduce the set
\begin{equation}\label{d3}
 \mathcal{X}=\{ v \in L_\infty(\T^3 \times T_{per}\S^1) : \div v=0 \mbox{ and } \|v\|_{L_\infty} \leq \delta\}.
\end{equation}
Let $\bar v \in \mathcal{X}$,
then we consider the linearization of (\ref{NSE})
\begin{equation}\label{NSE2}
 \begin{array}{l}
  v_t  -\nu \Delta v  + \nabla p= -\div(\bar v \otimes v) + F,\\[6pt]
\displaystyle \div v =0, \qquad \int_{\T^3} v(x,t)dx=0.
 \end{array}
\end{equation}
This process introduces a map $\mathcal{T}(\bar v)=v$. We will show $\mathcal{T}:\mathcal{X}\to \mathcal{X}$ and $\mathcal{T}$ is compact. As a result
this will establish existence of
a time periodic  solution to the nonlinear system (\ref{NSE}).

Similar to the 2d case we set
\begin{equation}\label{d4}
 v=V - \frac{1}{\Omega} H(x,\Omega t).
\end{equation}
In case (A) we  take $H=\sum_{k\in\Z^3} \frac{\hat{F}_k}{ik_1} e^{ikx}e^{i\Omega t k_1}$ and for case (B)  we take $H=-\cos \Omega t f(x)$, thus $\d_t \frac{1}{\Omega} H=F$.
Then we get
\begin{equation}\label{NSE3}
 \begin{array}{l}
\displaystyle  V_t  -\nu \Delta V  + \nabla p= -\div(\bar v \otimes V) + \div(\bar v \otimes \frac{1}{\Omega}H) +
  \frac{\nu}{\Omega} \Delta H.\\[6pt]
\div V =0,
 \end{array}
\end{equation}

The existence of solutions to (\ref{NSE3}) is sketched  in Appendix.
The estimates for solutions to (\ref{NSE3}) are done in  the domain $\T^3 \times \S^*$ with $\S^*=([\frac{1}{T_{per}}]+1)T_{per}\S^1$
just in order  to avoid the possible problem with smallness of $T_{per}$.
Provided we solved  system (\ref{NSE2}) with $\bar v\in \mathcal{X}$, we want to find a suitable estimate guaranteeing $v$ in $L_\infty$.
Here we work with the Slobodeckii spaces $W^{1,1/2}_p(\T^3 \times \S^*)$ \cite{A,LSU}.
In the Appendix we explain the details. Then we find the following inequality for  system (\ref{NSE3})
\begin{equation}\label{d5}
 \|V\|_{W^{1,1/2}_p(\T^3 \times \S^*)} \leq C_\nu \|\bar v V, \; \frac{1}{\Omega} \bar v H, \; \frac{\nu}{\Omega} \nabla H\|_{L_p(\T^3 \times \S^*)}.
\end{equation}
Now we use $\Omega$  large enough to guarantee the smallness of $\delta$  such that $C_\nu \|\bar v\|_{L_\infty} \leq 1/2$,
 (observe the norm $\|V\|_{W^{1,1/2}_p}$ contains $\|V\|_{L_p}$, as well), then
\begin{equation}\label{d6}
 \|V\|_{W^{1,1/2}_p(\T^3 \times \S^*)} \leq \frac{C_\nu (\delta  + \nu)}{\Omega} G.
\end{equation}
Next, we note that if $p>5$ then the space $W^{1,1/2}_p(\T^3 \times \S^*)$ is compactly imbedded in $L_\infty(\T^3 \times \S^*)$ \cite{BIN}, Chap XII.
 Consequently we have
\begin{equation}\label{d7}
 \|V\|_{L_\infty(\T^3 \times T_{per}\S^1)} \leq \frac{C(\delta  + \nu)}{\Omega} G.
\end{equation}
The constant is independent from $T_{per}$, since (\ref{d5}) is considered on $\T^3 \times \S^*$.
Therefore (\ref{d4}) implies that
\begin{equation}\label{d8}
 \|v\|_{ L_\infty(\T^3 \times T_{per}\S^1)} \leq C(\frac{(\delta  + \nu)}{\Omega} +\frac{1}{\Omega})G \leq \delta={C^*}\frac{1+\nu}{\Omega}G
\end{equation}
which is guaranteed for $\Omega$ larger than $\Omega_0(\nu,G)$.

Thus, we $\mathcal{T}:\mathcal{X}\to \mathcal{X}$ is compact. Hence the Schauder theorem implies existence of a fixed point of the map $\mathcal{T}$,
what yields existence of a time periodic solution to the original system (\ref{NSE}). Theorem 3 is proved.

\smallskip

\noindent
{\bf Attraction of weak solutions -- the 3d case}

\smallskip

{\it {\bf Theorem 4.} Let $v_0\in L_2(\T^3)$ be a divergence-free vector field, and let $v(t)$ be a Leray-Hopf weak solution to (\ref{NSE}) with initial datum $v_0$.
Then
\begin{equation}\label{f1}
 \|v(t)-v_{per}(t)\|_{L_2(\T^3)} \to 0, \mbox{ as } t \to \infty,
\end{equation}
where $v_{per}$ is the time periodic solution given by Theorem 3, provided $\Omega$ is large enough.
}

\smallskip

{\bf Proof.} Since weak solutions are not known whether they satisfy the energy equality, in the three-dimensional case, it will not be possible for us to follow the same arguments as
in the proof of Theorem 2.  However, since we are considering here   Leray-Hopf weak solutions, then, on one hand,  $v(t)$  satisfies following strong energy inequality
\begin{equation}\label{g1}
 \|v(t)\|^2_{L_2(\T^3)}+2\nu \int_s^t \|\nabla v(\tau)\|^2_{L_2(\T^3)} d \tau\leq \|v(s)\|^2_{L_2(\T^3)} + 2 \int_s^t (F,v)d\tau,
\end{equation}
for all $t>0$ and a.e. $s$ such that $t>s\geq 0$.
On the other hand, time periodic solutions are  regular solutions,  thus they do satisfy the energy   equality
\begin{equation}\label{g2}
 \|v_{per}(t)\|^2_{L_2(\T^3)}+2\nu \int_s^t \|\nabla v_{per}(\tau)\|^2_{L_2(\T^3)} d \tau= \|v_{per}(s)\|^2_{L_2(\T^3)} + 2 \int_s^t (F,v_{per})d\tau,
\end{equation}
for every $t \geq s \geq 0$.
To obtain an estimate for $\|v(t)-v_{per}(t)\|_{L_2(\T^3)}$ we use the observation that
\begin{multline}\label{g3}
 \|v(t)-v_{per}(t)\|_{L_2(\T^3)}^2=(v(t)-v_{per}(t),v(t)-v_{per}(t))=
\\
(v(t),v(t)) + (v_{per}(t),v_{per}(t)) - (v(t),v_{per}(t)) - (v_{per}(t),v(t)).
\end{multline}
To control the last terms we use the weak formulation for the solutions $v$ and $v_{per}$. Specifically, since $v_{per}$ is a sufficiently smooth we are allowed, on the one hand,  to use it as a test function in the weak  formulation for the weak solution $v$ to obtain
\begin{equation}\label{g4}
2\d_t(v,v_{per}) -2(v,v_{per,t}) + 2\nu (\nabla v,\nabla v_{per})
+ 2 (v\cdot \nabla v, v_{per}) = 2(F,v_{per}).
\end{equation}
On other hand, since $v_{per}$ is regular enough solution and the equation holds in $L_2(\T^3 \times T_{per}\S^1)$, we can multiply by $v$ and integrate over $\T^3$ to infer
\begin{equation}\label{g5}
 2(v_{per,t},v) + 2\nu (\nabla v_{per},\nabla v)
+ 2 (v_{per}\cdot \nabla v_{per}, v) = 2(F,v).
\end{equation}
Both (\ref{g4}) and (\ref{g5}) are meant in the distributional sense in time.
We add (\ref{g4}) and (\ref{g5}) and integrate over time interval $(s,t)$, and obtain
\begin{multline}\label{g6}
 2(v(t),v_{per}(t))+2\nu \int_s^t (\nabla v,\nabla v_{per})+(\nabla v_{per},\nabla v) d\tau \\
+ 2\int_s^t (v\cdot \nabla v, v_{per})+(v_{per}\cdot \nabla v_{per}, v)d\tau = \\
2(v(s),v_{per}(s))+2 \int_s^t (F,v + v_{per})d\tau
\end{multline}
Taking (\ref{g1})+(\ref{g2})-(\ref{g6}) we obtain
\begin{multline}\label{g7}
\|v(t)-v_{per}(t)\|_{L_2(\T^3)}^2 + 2\nu \int_s^t \|\nabla (v(t) - v_{per}(t))\|_{L_2(\T^3)}^2 d \tau \\
\leq \|v(s)-v_{per}(s)\|_{L_2(\T^3)}^2+ 2\int_s^t (v\cdot \nabla v, v_{per})+(v_{per}\cdot \nabla v_{per}, v)d\tau .
\end{multline}

Let $\delta v(t)= v(t) - v_{per}(t).$
 Consider the term
\begin{equation}\label{g8}
 (v\cdot \nabla v, v_{per})+(v_{per}\cdot \nabla v_{per}, v).
\end{equation}
We have

\begin{multline}\label{g9}
 (v\cdot \nabla v, v_{per})+(v_{per}\cdot \nabla v_{per}, v)=(v\cdot \nabla v, v_{per})-(v_{per}\cdot \nabla v, v_{per})=\\
(v\cdot \nabla v_{per}, v_{per})-(v_{per}\cdot \nabla v_{per}, v_{per})+(v\cdot \nabla \delta v, v_{per})-(v_{per}\cdot \nabla \delta v, v_{per})=\\
0+0 + (\delta v \cdot \nabla \delta v, v_{per}).
\end{multline}
Two first terms vanished. Next we note that
\begin{equation}\label{g10}
 |\int_{\T^3} \delta v \cdot \nabla \delta v \,  v_{per} dx|\leq \frac{\nu}{2} \|\nabla \delta v\|^2_{L_2(\T^3)} +
\frac{C}{\nu} \|v_{per}\|^2_{L_\infty} \|\delta v\|^2_{L_2(\T^3)}.
\end{equation}
Next, the Poincar\'e inequality yields $\|\delta u\|_{L_2(\T^3)} \leq \|\nabla \delta u\|_{L_2(\T^3)}$. Moreover, we observe  that  $\|v_{per}\|_{L_\infty}\leq  \delta$ -- see the proof of Theorem 3,
 with $\delta$ very small so that $\frac{C}{\nu} \delta^2 \leq \frac{\nu}{4}$ which holds as $\Omega$ is sufficiently large. Using the above to finally
obtain
\begin{equation}\label{g11}
 \|\delta v(t)\|_{L_2(\T^3)}^2 + \nu \int_s^t \|\nabla \delta v(t)\|_{L_2(\T^3)}^2 d \tau \\
\leq \|\delta v(s)\|_{L_2(\T^3)}.
\end{equation}
Again using the Poincar\'e inequality we obtain
\begin{equation}\label{g12}
 \|\delta v(t)\|_{L_2(\T^3)}^2 +  \nu \int_s^t \|\delta v(t)\|_{L_2(\T^3)}^2 d \tau \\
\leq \|\delta v(s)\|_{L_2}^2,
\end{equation}
for all $t>0$ and a.e. $s\geq 0$. Simple analysis of (\ref{g12}) implies
\begin{equation}\label{g13}
 \|\delta v(t)\|_{L_2(\T^3)}\leq Ce^{-c \nu \, t}.
\end{equation}
Theorem 4 is proved. As a corollary we obtain the fact that time periodic  established by Theorem 3 are unique for sufficiently large $\Omega$.

\section{Numerical results for 2D case}

The numerical results presented in this section focus on  the system \eqref{EeS},
i.e.
\begin{subequations}
\begin{gather}
  \label{eq:2dns}
  v_t + v\cdot\nabla v - \epsilon\Delta v + \alpha(v-\Omega\hat{e}_2)  +\nabla p = F(z),\\
  {\rm\ div\,}{v} = 0,\\
 \frac{1}{|\T^2|} \int_{\mathbb{T}^2}{v_0\,dz} = \Omega\hat{e}_2,
\end{gather}
\end{subequations}
where $v\colon[0,\infty)\times\mathbb{T}^2\to\mathbb{R}^2$, $z=(x,y)$.
In the following sections we are concerned with the
numerical investigation of the dependence of the long-time qualitative behavior of the solutions of the above system on the parameter $\Omega$, for a given particular external forcing term $F(z)$. In the sequel we are going to assume \eqref{i2},
i.e., that average of the forcing term  over $\mathbb{T}^2$ is zero. In turn, this implies that the spatial average of the solutions remain constant.

In view of the theoretical analysis, presented in the previous sections, it follows that when  the values  of $\Omega$ exceed  certain critical value
implies the stabilization of the evolutionary problem \eqref{eq:2dns}.
More precisely, for a given particular forcing $F(z)$  the  stationary problem
\eqref{eq:2dns} might have multiple solutions, however, after increasing $\Omega$ beyond certain critical value one obtains a unique stationary solution.

\subsection{Numerical investigation particular setting}

In our numerical investigation we focus on a particular case study of the Kolmogorov flow  that was discussed in details at the end of section~2. Specifically we consider system \eqref{eq:omegaExpOld}  in the flat torus $\mathbb{T}^2_\beta$. As we have discussed earlier, in the end of section~2, the global stability of time period solutions of \eqref{eq:omegaExpOld} is equivalent to the global stability of stationary solutions of \eqref{vortSystem}. For this reason we focus in the next section the study of the bifurcation diagram of stationary solutions of \eqref{vortSystem}.

\subsection{Stationary problem bifurcation analysis}
After dropping the tilde system \eqref{vortSystem} is given by
\begin{gather}
\label{eq:rotBefore}
\omega_t + (v + \beta\Omega\hat{e}_2)\cdot\nabla\omega-\epsilon\Delta\omega+\alpha\omega = \lambda\cos(\frac{y}{\beta}),
\end{gather}
subject to periodic boundary condition, with basic domain $\mathbb{T}_\beta^2=[0,2\pi]\times [0,2\pi\beta]$.
In this section we present our numerical investigations of the stationary problem of \eqref{eq:rotBefore}:
\begin{gather}
\label{eq:omega}
  (v + \beta\Omega\hat{e}_2)\cdot\nabla\omega-\epsilon\Delta\omega+\alpha\omega = \lambda\cos(\frac{y}{\beta}),
\end{gather}
subject to periodic boundary condition, with basic domain $\mathbb{T}_\beta^2=[0,2\pi]\times [0,2\pi\beta]$.

\begin{figure}[htbp]
  \centering
  \includegraphics[width=0.75\textwidth]{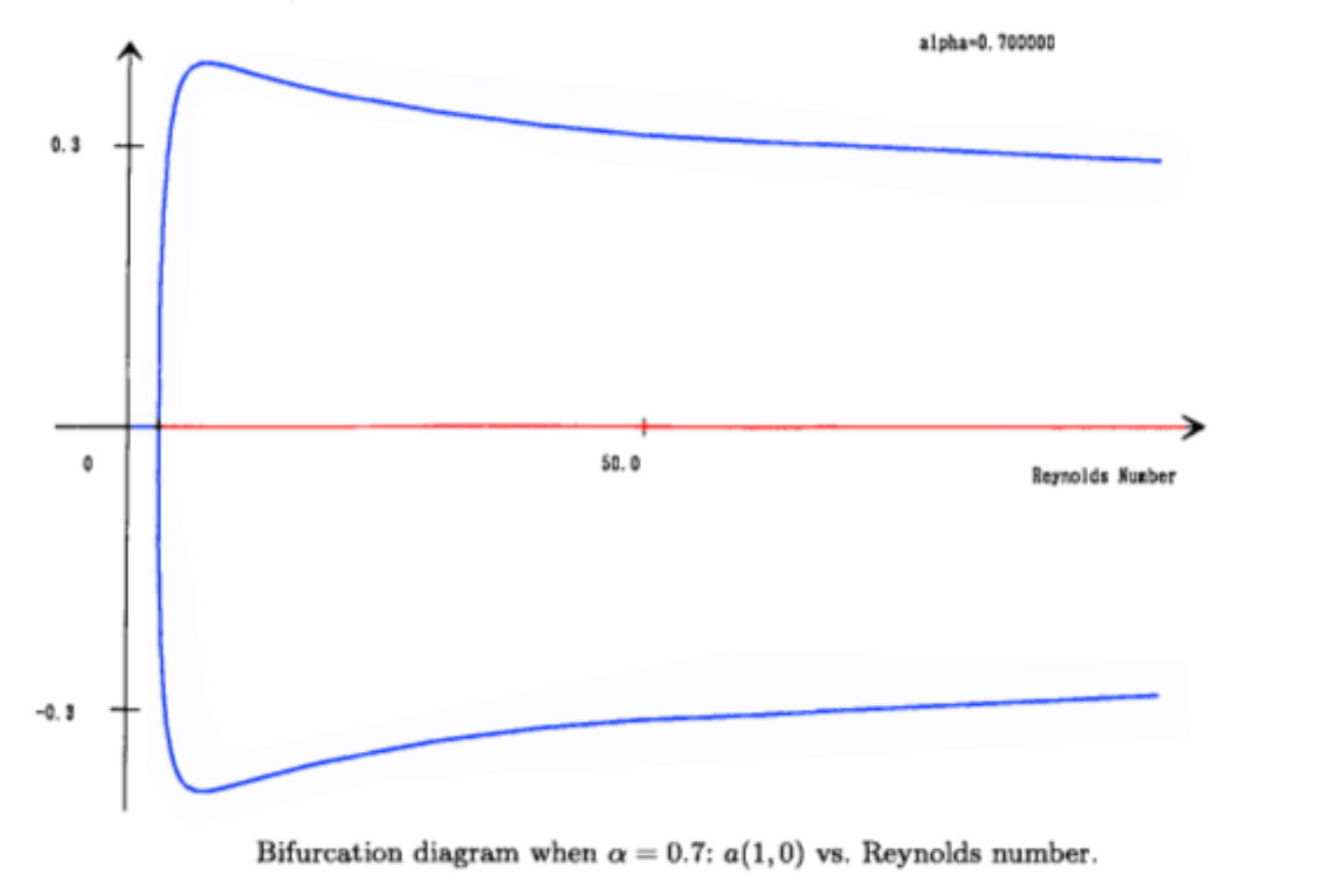}
  \caption{The bifurcation diagram for the problem \eqref{eq:rotBefore} with
    $\beta=0.7, \Omega=0, \epsilon=1, \alpha =0$ (reproduced from \cite{OS}), the scaling parameter here is the Reynolds number
    defined by the authors
    $Re \simeq \frac{\lambda}{\epsilon^2\beta^3}$. Let $a(1,0)$ denotes the $\omega$'s Fourier coefficient corresponding to
    $\exp{ix}$ basis function. The stable solution (in blue) having $a(1,0)=0$ becomes unstable (in red) at particular
    critical value of the Reynolds number, where two new stationary solutions are being born.}
  \label{fig:bif}
\end{figure}

In the case $\beta = 1$ (the square) and when $\Omega=0$
system \eqref{eq:rotBefore}  admits a globally stable stationary solution (called trivial solution) \cite{CFT,M}. Consequently this globally stable stationary  solution does not undergo any bifurcation when $\lambda$ is increased.
On the other hand, numerical experiments in \cite{OS} show that for $\beta\in(0,1)$  system \eqref{eq:omega}, when
$\Omega=0$, the trivial stationary solution undergoes a pitchfork bifurcation (see also \cite{BV}). In the remaining part of this section we restrict our attention to the
particular case $\beta=0.7$, for which we present bifurcation diagram on Figure~\ref{fig:bif} (reproduced from \cite{OS}).

Looking at Figure~\ref{fig:bif} it is evident that  problem \eqref{eq:omega}, with $\Omega=0$, exhibits unique stationary
solution for $\lambda$ values smaller than a critical value (we denote it by $\lambda_0$) -- the point of the pitchfork
bifurcation, at which two branches of stable solutions are born. Let $\lambda_0(\Omega)$ denote the point of
the pitchfork bifurcation in problem \eqref{eq:omega}, and let $\omega(\lambda_0(\Omega))$ denote the solution at
bifurcation point.

\begin{figure}[htbp]
  \centering
        \begin{subfigure}[b]{0.7\textwidth}
          \includegraphics[width=\textwidth]{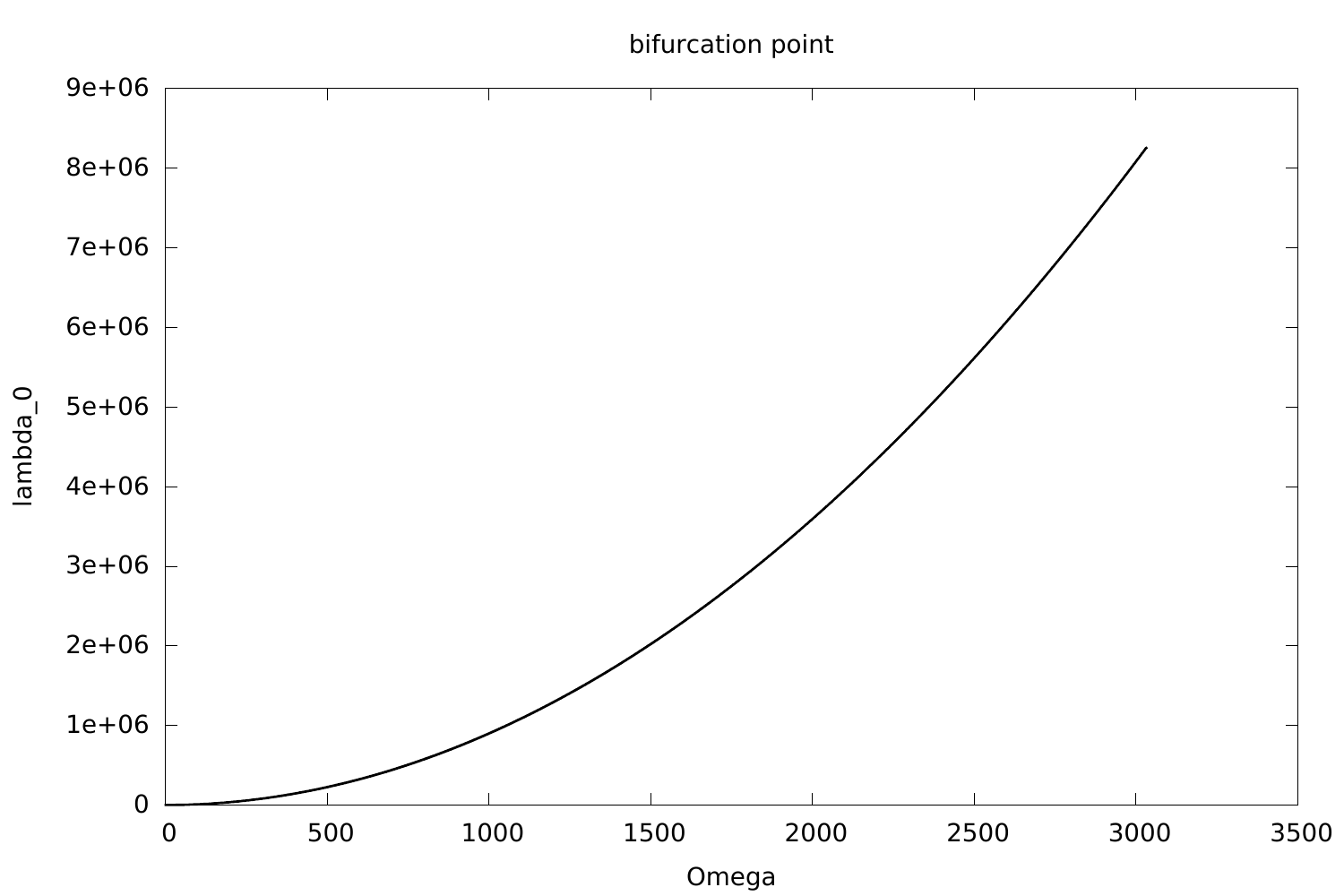}
          \caption{diagram showing the bifurcation point $\lambda_0$ with respect to $\Omega$. The presented graph is approximately
          $
            \lambda_0(\Omega)=0.89851\cdot\Omega^2 + 0.0014294\cdot\Omega + 1,
          $
          obtained by least squares fitting with root mean square of residuals equal to $0.612266$}
        \label{fig:bifurcation}
        \end{subfigure}
        \begin{subfigure}[b]{0.7\textwidth}
          \includegraphics[width=\textwidth]{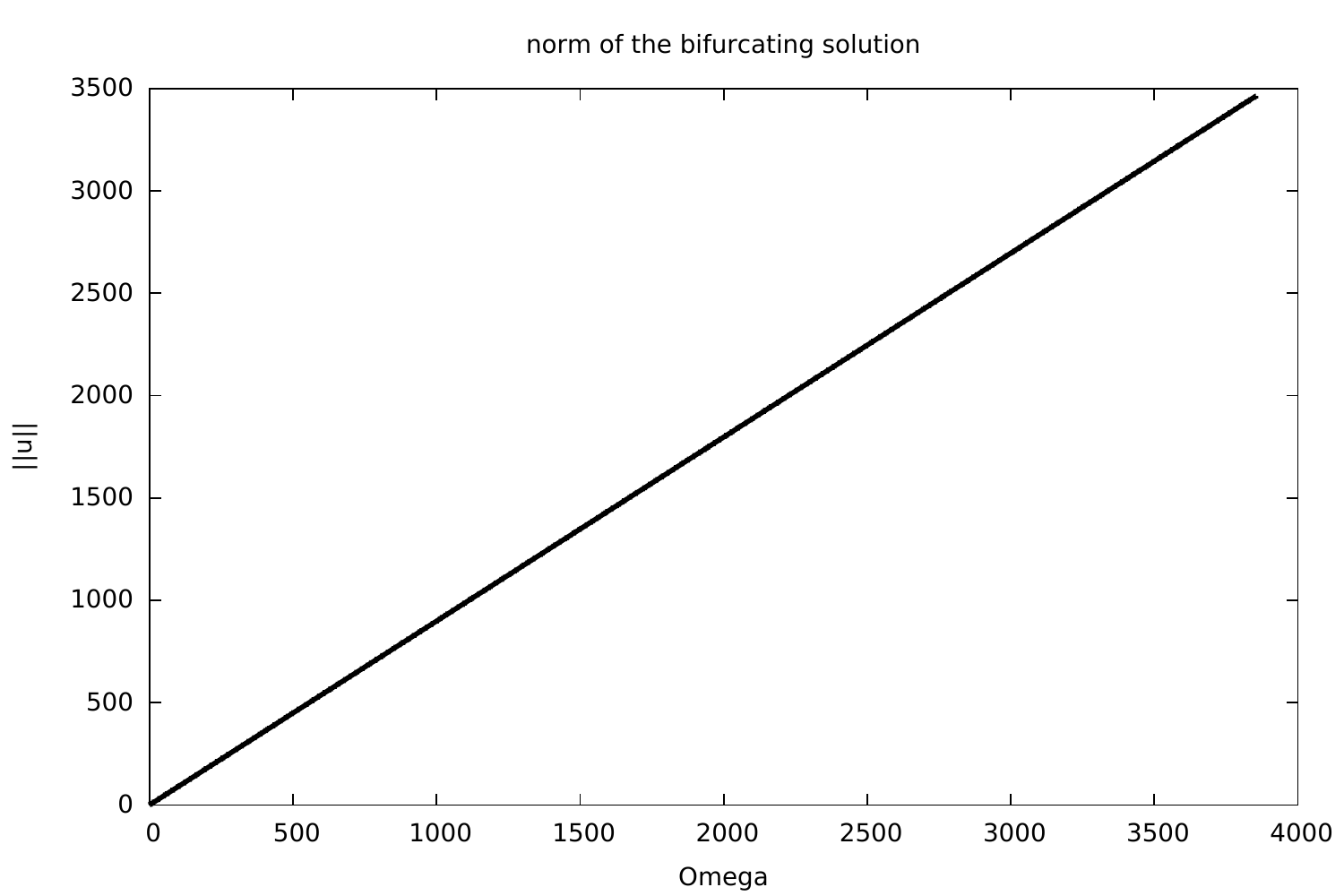}
          \caption{diagram showing $\|\omega(\lambda_0)\|$ ($L_2$ norm of the solution at bifurcating point) with respect to $\Omega$ }
        \end{subfigure}
        \caption{Bifurcation diagrams for the problem \eqref{eq:rotBefore} with $\alpha=0$, and $\epsilon=1$.}
       \label{fig:nu}
\end{figure}

We investigate here the dependence of $\lambda_0(\Omega)$, and the dependence of $\|\omega(\lambda_0(\Omega))\|$ --
the solution at bifurcation point $L_2$ norm on $\Omega$.
The numerical tests are in agreement with the theory presented in the theoretical  part of this paper, from which it follows
that the region  of the parameter $\lambda$, for which one has unique stable stationary solution of \eqref{eq:omega}, is enlarged while $\Omega$ increases,
and that $\|\omega(\lambda_0(\Omega))\|$ should increase with the order of magnitude lower than that of $\lambda_0$.
Figure~\ref{fig:nu} agrees with theoretical derivations of \eqref{bb5} for equation \eqref{eq:omega} showing $\|\omega\|_{L_2}\sim \lambda \Omega^{-1}$
($G$ from \eqref{bb5} is proportional to $\lambda$). Presented numerical results indicate also that the lost of uniqueness occur for
$\lambda \gtrsim \Omega^2$ hence the norm of the solution at bifurcation point is of linear growth in $\Omega$.

In Figure~\ref{fig:nu} we present the calculated results for the problem \eqref{eq:omega}  with $(\epsilon,\alpha)=(1,0)$, we
skip here results for other choices of $(\epsilon,\alpha)$, as we did not observe any qualitative difference
in this case.

\subsection{Effect of stabilization}
This part illustrates the stabilization effect for problem \eqref{eq:rotBefore}. For particular initial conditions provided later on
we integrate in time the evolution equations \eqref{eq:rotBefore}. In \eqref{eq:rotBefore} we force the second mode
(the forcing is $\lambda\cos(\frac{2y}{\beta})$), as we observe a rich dynamics for that case.
For a fixed $\lambda$ we compare $\Omega=0$ case with $\Omega$ large.
As a result we obtain the stabilization effect with exponential convergence rate for the latter case, as Theorem~2 predicts.

In order to numerically integrate \eqref{eq:rotBefore} forward in time we invoke standard numerical integrator.
We write $\omega$ in (complex) Fourier basis, i.e., $\omega(t,z)=\sum{a_k(t)\exp{i(k,z)}}$. We consider a Galerkin
approximation of the infinite system of ODEs including only modes $a_k$ with $k$ such that $|k|_\infty\leq N$. We call
$N$ the approximation dimension. In the presented experiment we fixe $N=13$, this choice is motivated by the fact that this approximation
dimension represents well the dynamics of the PDE, we validate this by checking that for a larger dimension ($N=19$) the obtained
results are qualitatively the same (not provided here). To perform
time integration procedure we use the Taylor method, the time step is selected adaptively, is maximized under constraint such that the local error do not exceed the machine precision.
In our actual computations we fix Taylor's method order to $15$, which is relatively high order as for a high dimensional system, however,
in our case it provides an efficient procedure.

We describe the following numerical experiment.

We fix $\lambda=100$, $\beta=0.75$, $N=13$, order of the Taylor method is $15$. 
We pick four initial conditions and integrate the equation on the time interval $[0,5]$.
\begin{enumerate}
  \item Initial condition I -- $\omega^I_0(x,y) = 60\cos(x+y)$; it is attracted by a periodic orbit of $L_2$ norm in $(19.7, 20)$,

  \item Initial condition II -- $\omega^{II}_0(x,y) = 2\cos(x)$; it is attracted by a stationary solution of $L_2$ norm approximately $9.68043$,

  \item Initial condition III -- $\omega^{III}_0(x,y) = 2\cos(2x)$; it is attracted by a stationary solution of $L_2$ norm approximately $14.2384$,

  \item Initial condition IV -- $\omega^{IV}_0(x,y) = 0$; it is attracted by a stationary solution of $L_2$ norm approximately $25$.

\end{enumerate}

\begin{figure}[htbp]
  \centering
        \begin{subfigure}[b]{0.75\textwidth}
          \includegraphics[width=\textwidth]{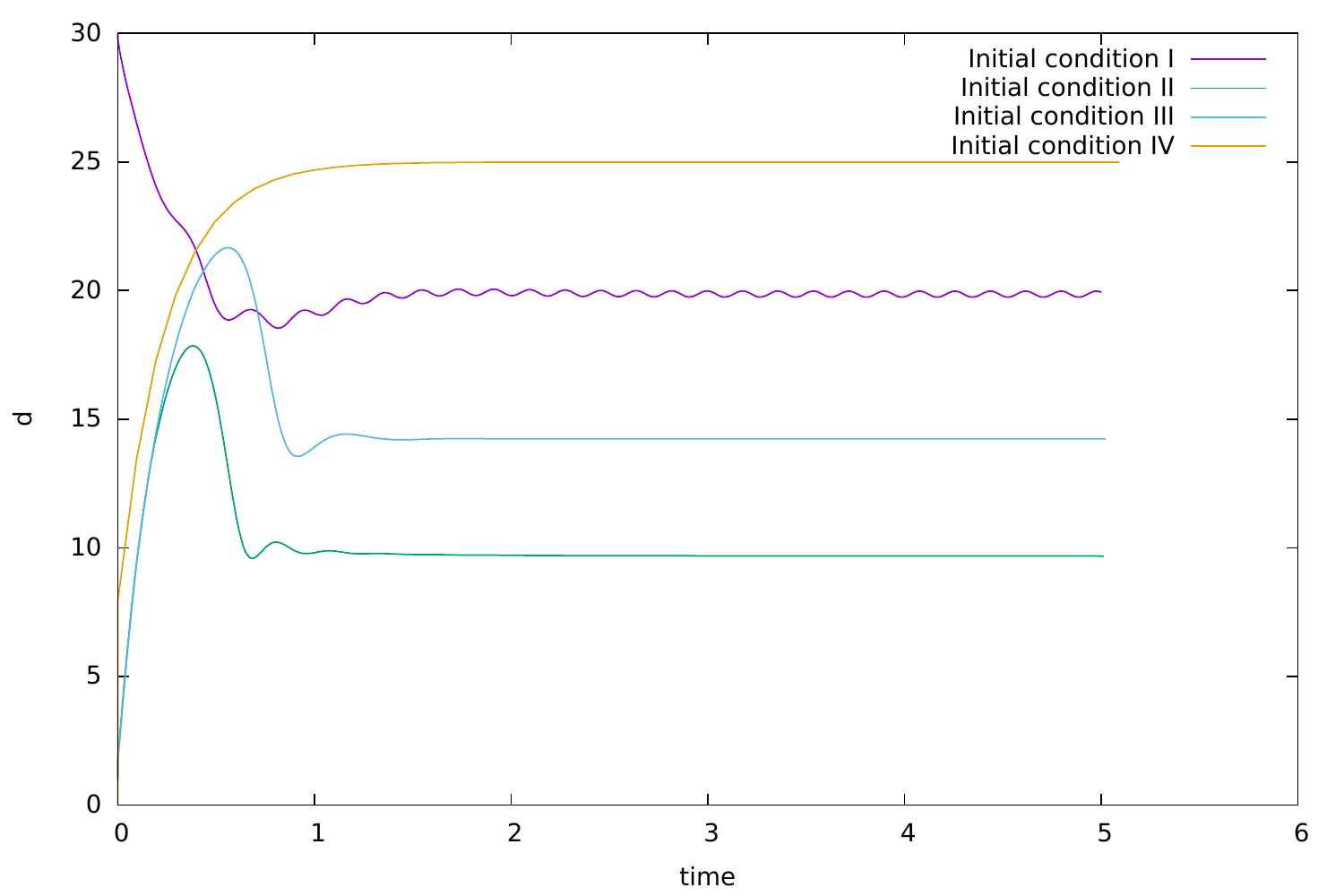}
          \caption{Case of $\Omega = 0$}
          \label{a}
        \end{subfigure}
        \begin{subfigure}[b]{0.75\textwidth}
          \includegraphics[width=\textwidth]{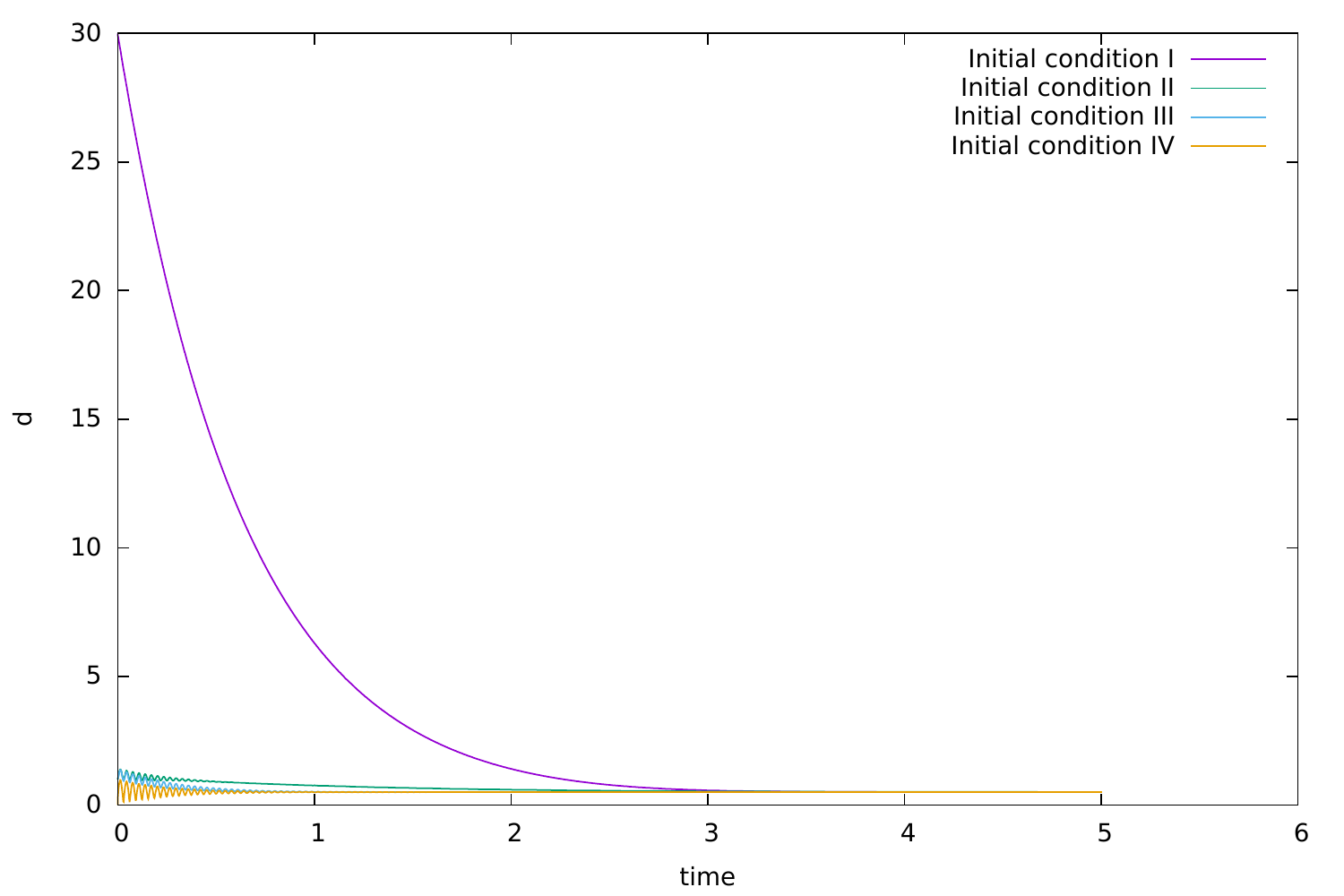}
          \caption{Case of $\Omega = 100$}
          \label{b}
        \end{subfigure}
        \caption{Results from integrating in time the equation \eqref{eq:rotBefore} with $\Omega=0$, and $\Omega=100$,
using four different initial conditions $\{\omega^I_0$, $\omega^{II}_0$, $\omega^{III}_0$, $\omega^{IV}_0\}$.
The $L_2$ norm of the solution with respect to time is plotted, we denote $d=\|\omega(t)\|_{L_2}$.}
       \label{fig:omegas}
\end{figure}

Observe that for $\Omega=100$ all of the initial conditions, even the periodic orbit case, are being attracted by the same stationary solution,
so the stabilization is achieved for $\Omega \gtrsim \sqrt{\lambda}$, which agrees with our expectations.

\subsection{Evolutionary problem convergence rate analysis}
In this section we present the results of our investigation of the global convergence to the unique stationary solution, for large initial values,  of the general time dependent solutions of the evolution equation \eqref{eq:rotBefore}.

First we fix the forcing amplitude $\lambda=1$, in \eqref{eq:rotBefore}, and we consider the  initial value $\omega_0=\Omega(\sin{x}+\cos{x})$,
 $\Omega$ here is  the large parameter and is the amplitude of the initial value.
We recorded the time needed for $\omega_0$ to be attracted by the stationary solution
of \eqref{eq:rotBefore}, with $\lambda=1$. We stopped our numerical integration procedure at time $T$
when
\begin{equation*}
  \|v(T)\cdot\nabla\omega(T)-\epsilon\Delta\omega(T)+\alpha\omega(T)-\lambda\cos(\frac{y}{\beta})\|_{L_2}\leq 10^{-5}.
\end{equation*}

\begin{figure}[htbp]
  \centering
        \begin{subfigure}[b]{0.5\textwidth}
          \includegraphics[width=\textwidth]{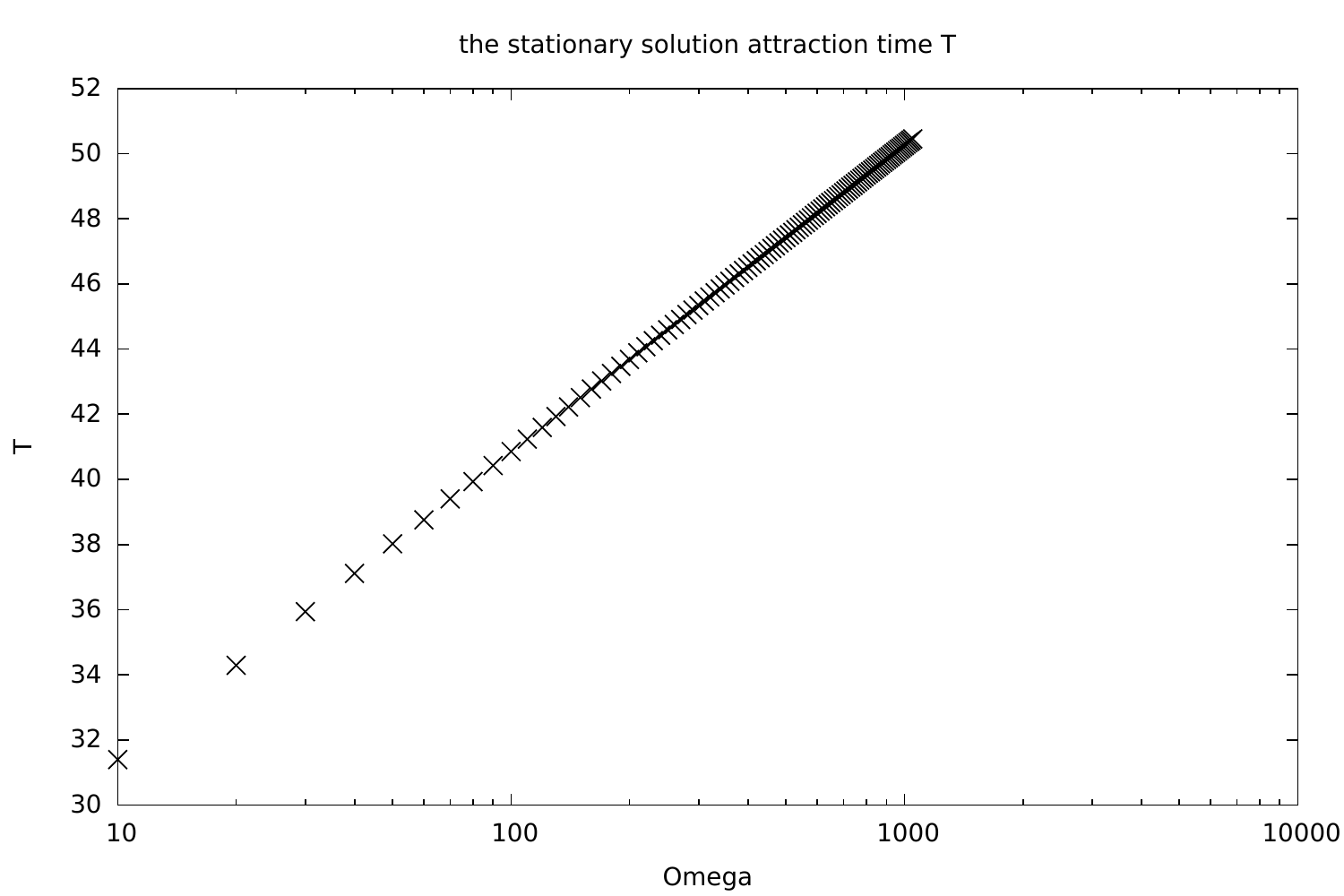}
\caption{Logscale plot of the stationary solution attraction time of (denoted $T$) with respect to
$\Omega$ ($\lambda=1$ is fixed). The initial condition for a fixed $\Omega$ is $\omega_0=\Omega(\sin{x}+\cos{x})$.
Apparent linear growth of $T$ matches the exponential convergence established in Theorem~2.}
\label{fig:convergence}
        \end{subfigure}
        \begin{subfigure}[b]{0.45\textwidth}
          \includegraphics[width=\textwidth]{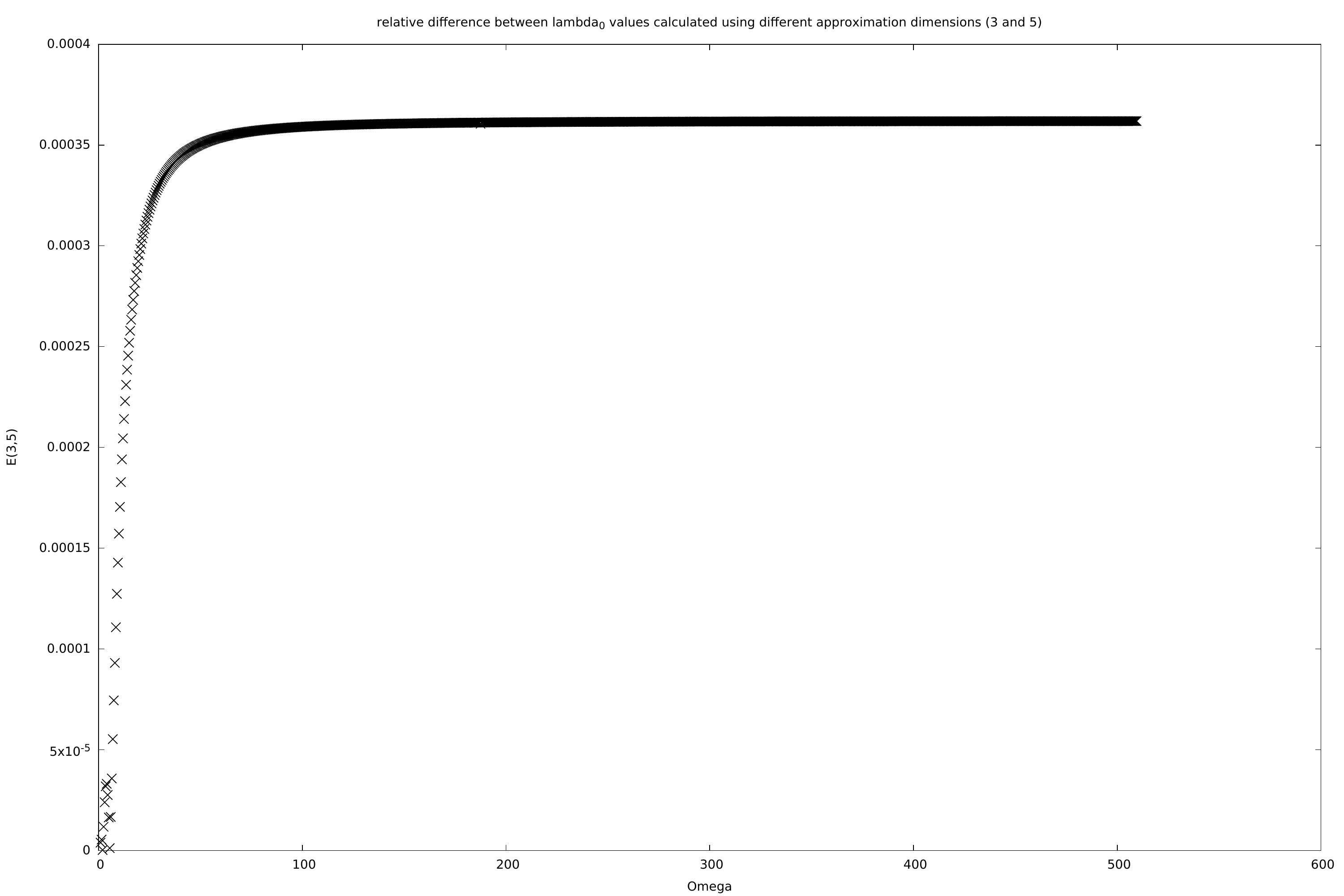}
\caption{Diagram presenting\\$E(\Omega,3, 5)=\left| \left(\lambda_0(\Omega,3)-\lambda_0(\Omega,5)\right) / \lambda_0(\Omega,3) \right|$, the relative difference between $\lambda_0(\Omega,3)$, and $\lambda_0(\Omega,5)$, the parameter here is $\Omega$.\\ \ \\}
\label{fig:error}
        \end{subfigure}
\end{figure}

 In Figure~\ref{fig:convergence} we present the results
for the case $(\epsilon,\alpha)=(1,0)$ (other cases were qualitatively very similar).
Figure~\ref{fig:convergence} is plotted in the logscale, and the apparent linear growth of $T$
matches the exponential convergence established in Theorem~2.

To perform the numerical time integration we invoke the same techniques as in the previous section,
however, for this particular choice of initial condition and $\lambda=1$ the dynamics is apparently low dimensional,
so we fix $N=3$.
This low Galerkin approximation may seem not sufficient.
Therefore to argue that the results for larger Galerkin approximation, in this case, do not differ qualitatively from the obtained results,
we perform an additional numerical test in which we measure the relative difference between two approximation of
dimensions $N$ and $M$ by
  $E(\Omega,N,M) = \left| \left(\lambda_0(\Omega, N)-\lambda_0(\Omega, M)\right) / \lambda_0(\Omega, N) \right|$,
where $\lambda_0(\Omega, N)$ is the bifurcation point for particular $\Omega$,
calculated using the approximation dimension $N$.
We present the obtained diagrams for $N=3$, and $M=5$ in Figure~\ref{fig:error}. Clearly, the values shown remain
essentially constant for larger $\Omega$ values, and does not exceed $0.0004$, which supports our claim that
it is enough to use a small approximation dimension to provide a qualitatively correct illustration.


\subsection{Conclusions from numerical experiments and future work}
The goal of this section was to present a numerical investigations of a simple case, as an illustration of the theoretical results presented in the theoretical sections of this paper.
Obviously, the numerical results match the theoretical predictions. All theorems in this paper are about periodic solutions, but in order to obtain the numerical results we always reduce the problem to the stationary case. This is imposed by the fact that for high values of $\Omega$ the periodic solutions oscillate rapidly, which is a major obstacle for
 numerical integration in time. Due to the equivalence between the two problems, as we have indicated in this section, the conclusions from the numerical results are meaningful for the
case of oscillating rapidly periodic solutions, although the computations are performed for the stationary case. There exist several numerical methods which probably allow to treat the case with
rapid oscillations directly, but our current goal was solely to provide an illustration for the theoretical results  established in this paper, rather than invoking sophisticated numerical methods
to deal with rapid oscillations generated by large values of $\Omega$. We leave the task for future  research.

\section{Appendix}

In this part we explain the construction of time periodic solutions for the ``linearized" problem. We establish this for the three-dimensional  case, which is an essential step in the proof of Theorem 3. The case for Theorem 1 is almost the same.

Having $\bar v \in \mathcal{X}$ we consider (\ref{NSE3}) in the following form
\begin{equation}\label{x1}
 \begin{array}{l}
  V_t- \nu \Delta V +\nabla p = -\div (\bar v \otimes V) + \div (\bar V \otimes \frac{1}{\Omega} H) + \frac{\nu}{\Omega} \Delta H,\\
\div V=0.
 \end{array}
\end{equation}
The issue of existence for the above system lays in the classical theory. The simplest approach is through Fourier methods using series in time and space to the
linear system with given right-hand side. We acts on the domain $\T^3 \times T_{per} \S^1$ and we represent the solution in  the form
\begin{equation}\label{x2}
 V^{(j)}(x,t) \sim \sum_{l \in \Z, k \in \Z^3} V_{lk}^{(j)} e^{ilt T_{per}} e^{ikx}, \qquad j=1,2,3.
\end{equation}
The solvability of the system for the finite dimensional approximation is clear, so we need just a good estimate which allows to pass to the limit.
But the energy estimate is allowed to be used in the chosen framework, so we get
\begin{equation}\label{x3}
 V \in L_2(T_{per} \S^1;H^1(\T^3)),
\end{equation}
with the a priori estimate
\begin{equation}\label{x4}
 \|V\|_{L_2(T_{per}\S^1;H^1(\T^3))} \leq C (\frac{1}{\Omega} \|\bar v\|_{L_\infty} \|H\|_{L_2(\T^3\times T_{per}\S^1)} +
 \frac{\nu}{\Omega} \|\nabla H\|_{L_2(\T^3\times T_{per}\S^1)}).
\end{equation}
The construction by approximation in the time periodic functions ensures the solution $V$ is $T_{per}$ - periodic.
The form of (\ref{x2}) guarantees the periodicity in time and space. The term $\div (\bar v \otimes V)$
can be treated as a perturbation, and thanks to the estimate (\ref{x4}) we obtain the existence to  system (\ref{x1}), too.

Next, we  improve the regularity of solutions $V$. Here we apply the maximal regularity result for the Stokes operator in the
 Slobodeckii spaces of type $W^{1,1/2}_p(\T^3 \times \S^*)$  \cite{A,BIN,Tr}.
For solutions to
\begin{equation}
 \begin{array}{l}
  V_t- \nu \Delta V +\nabla p = \div F,\\
\div V=0
 \end{array}
\end{equation}
the following estimate holds
\begin{equation}
 \|V\|_{W^{1,1/2}_p(\T^3 \times \S^*)} \leq C_\nu \|F\|_{L_p(\T^3 \times \S^*)}.
\end{equation}
The definition of the Slobodeckii space \cite{Sol,Tr} is by the following norm
\begin{equation}
 \|V\|_{W^{1,1/2}_p(\T^3 \times \S^*)}^p = \|V\|_{L_p(\T^3 \times \S^*)}^p+ \|\nabla V\|_{L_p(\T^3 \times \S^*)}^p+
\int_\Omega \int_{\S^*} \int_{\S^*} \frac{|V(x,t)-V(x,t')|^p}{|t-t'|^{1 + \frac{1}{2}p} } dtdt'dx.
\end{equation}
Thus we justify estimate (\ref{d5}).

\bigskip

\noindent{\bf Acknowledgments.}
The presented work has been done while J.C. held a post-doctoral position at Warsaw Center
of Mathematics and Computer Science, and more recently at Rutgers -- The State University of New Jersey, his research has been partly supported by Polish National Science
Centre grant 2011/03B/ST1/04780. The second author (P.B.M.) has been partly supported by National  Science  Centre  grant 2014/14/M/ST1/00108 (Harmonia).
The work of  E.S.T.  is supported in part by the ONR grant N00014-15-1-2333 and the NSF grants DMS-1109640 and DMS-1109645.

\end{document}